\documentclass[12pt]{article}

\usepackage[margin=1in,tmargin=1in]{geometry}

\usepackage{amsmath, amssymb, amsthm}
\usepackage{cancel}
\usepackage{amsfonts}
\usepackage{graphicx}
\usepackage[frame,ps,matrix,arrow,curve,rotate]{xy}
\usepackage{enumerate}
\usepackage{tikz}
\usepackage[colorlinks=true,citecolor=blue,linkcolor=blue,urlcolor=blue]{hyperref}
\usepackage{siunitx}
\usepackage{enumitem}
\usepackage{xspace}
\usepackage{centernot}
\usepackage{blkarray}
\usepackage{tikz-cd}
\usepackage{bm}

\usepackage{alphabeta}

\newtheorem{theorem}{Theorem}[section]

\newtheorem{corollary}[theorem]{Corollary}

\theoremstyle{definition}
\newtheorem{definition}[theorem]{Definition}

\theoremstyle{remark}

\newcommand{\mc}[1]{\mathcal{#1}}

\newcommand{\abs}[1]{\lvert#1\rvert}

\newcommand{\set}[1]{\{#1\}}

\renewcommand{\geq}{\geqslant}

\renewcommand{\leq}{\leqslant}
\renewcommand{\emptyset}{\varnothing}

\newcommand{\es}{\emptyset}

\DeclareMathOperator{\id}{id}

\DeclareMathOperator{\Spec}{Spec}

\newcommand{\bR}{\mathbb{R}}

\newcommand{\bQ}{\mathbb{Q}}

\newcommand{\bN}{\mathbb{N}}

\newcommand{\bZ}{\mathbb{Z}}

\newcommand{\cB}{\mathcal{B}}
\newcommand{\cN}{\mathcal{N}}
\newcommand{\cM}{\mathcal{M}}

\newcommand{\bSigma}{\bm{\Sigma}}
\newcommand{\bPi}{\bm{\Pi}}

\usepackage{blindtext}
\usepackage{titlesec}

\begin{document}
    \title{The set-theorist's real numbers}
    \author{Dhruv Kulshreshtha and Jamie Tappenden} 
    \maketitle

\begin{abstract} 
Descriptive set theorists often refer to elements of the Baire space, $\cN = {}^\omega\omega$, as ``reals.'' In this expository article, we present mathematical justification for this practice and highlight the advantages of working with $\cN$ over $\bR$ in the context of descriptive set theory.
\end{abstract}

\section{Introduction}

Descriptive set theory (DST), often described as the study of \textit{definable} (in this respect, ``well-behaved'') sets of real numbers, naturally interacts with various areas of non-foundational mathematics by examining the \textit{descriptive complexity} of sets arising in those fields.\footnote{We will explain more precisely the relevant sense of descriptive complexity in Section \ref{sec: favorable properties of Baire space}. Here, we are really referring to the field of \textit{classical DST}, whose roots trace back to the early twentieth-century works of Baire, Borel, Lebesgue, Lusin, Suslin, Sierpi\'nski, and others. In the 1950s, the classical theory was introduced to the machinery of computability theory---developed by Church, G\"odel, Kleene, Post, Turing, and others---by Addison, Mostowski, and others, leading to the development of the field of \textit{effective DST}. Since then, the field has expanded in several directions, from the study of \textit{definable} equivalence relations in \textit{invariant DST} to recent connections between \textit{descriptive combinatorics} and \textit{distributed algorithms}. Standard sources for classical, effective, and invariant DST include \cite{Kechris-classical-DST-1995,Moschovakis-DST,Gao-2009}, while an accessible introduction to the aforementioned recent connections can be found in \cite{Pikhurko-Borel-Combinatorics-2021,Bernshteyn2022}.} However, it is common for descriptive set theorists, when referring to ``real numbers''---often simply called ``reals''---to actually mean elements of Baire space, $\cN = {}^\omega\omega$, rather than elements of the set $\bR$ of ordinary real numbers.\footnote{See, for example, \cite[p.339]{Kechris-1973} and \cite[p.445]{Harrington-Kechris-1975}: ``Let $\omega = \set{0,1,2\ldots}$ be the set of natural numbers and $\bR = {}^{\omega}\omega$ the set of all functions from $\omega$ to $\omega$ or (for simplicity) \textit{reals}.'' Similar identifications are made in \cite[p.391]{Kechris-Moschovakis-1972}, \cite[p.405]{Martin-Kechris-1980}, \cite[p.147]{Steel-1981}, and \cite[p.379]{Kechris-Woodin-2008}. See also the definition of $\Gamma$-characteristic in \cite[p.3]{Blass-1993}, where Blass uses ``the descriptive set-theorist's usual version ${}^\omega\omega$ of the reals.''\label{footnote: examples where Baire space is called reals}
Blass also notes this identification in \cite[p.397]{Blass-cardinal-char-2010}, where he briefly comments on the meaning of ``continuum''. Blass also offers similar remarks about the identification of ``reals'' with elements of the Cantor space $\mc{C}$ (described in footnote \ref{footnote: Cantor space as reals}) in \cite[p.397]{Blass-cardinal-char-2010}.} 

Why is this adoption of $\cN$ as a representative for the real numbers a reasonable thing to do? The answer draws on assumptions of the sort that mathematicians call ``folklore'',  typically unstated and unexplained in print. As Moschovakis puts it: ``results in the theory are often proved just for $\cN$, with the (suitable) generalizations to other spaces and the reals in particular left for the reader or simply stated without proof'' \cite[p.4]{Moschovakis-DST}. Not only is it often unstated, but on the occasions where parts of it {\it{are}} stated, others are not. (In this village, the siblings in the story meet a mermaid but not trolls, in that village trolls but no mermaid,\ldots). Our purpose in this note is to serve as the Brothers Grimm for this particular folk tradition, to assemble scattered parts and flesh them out into a connected tale.

To this end, we explore the following questions: Why, in the context of DST, might one choose to code the reals as elements of Baire space? And what are the broadly understood but rarely stated principles justifying this practice?

This article is organized as follows: In Section~\ref{sec: preliminaries}, we review the relevant topological prerequisites. In Section~\ref{sec: favorable properties of Baire space}, we describe favorable properties of $\cN$ that may motivate set theorists to work with it in place of $\bR$. More specifically, we show how $\cN$ differs from $\bR$ in ways which, though relevant to characterizing $\bR$ in certain research contexts, prove fruitful for the directions that concern us. In Section~\ref{sec: Reducing R to Baire space} (respectively, Section~\ref{sec: reducing general case to Baire space}), we describe the ways in which $\cN$ is similar to $\bR$ (respectively, general Polish and quasi-Polish spaces), allowing us to reduce to the case of $\cN$.\medskip

\noindent\textbf{Acknowledgements.} We thank Andreas Blass for his insights throughout this project, especially for a conversation with the first author that substantially shaped Section~\ref{sec: Reducing R to Baire space}, and for his comments on an earlier draft. In addition, we thank Carolin Antos for helpful discussions and the writings of Paul Halmos for their influence on this project. The second author would like to thank John Burgess for initially introducing him to the basics of this folklore many years ago. The second author would also like to acknowledge a debt to his PhD supervisor Saul Kripke, who taught him much about descriptive set theory (ingeniously disguised as generalized recursion theory). The first author thanks Uri Andrews and Joseph S. Miller for their helpful comments and questions during his talk on this subject at the UW–Madison graduate logic seminar. We thank Luca Motto Ros for a helpful clarification concerning our discussion of the lexicographic order topology on $\cN$ in an earlier draft. We also thank Christian Rosendal for the reference \cite{Miller-Graph-Theory-DST-2012} to Benjamin Miller's proof of the Kechris--Solecki--Todorcevic $G_0$ dichotomy theorem, and Reed Solomon for the reference \cite{Marker-DST-notes} to David Marker's notes, where the proof of Theorem~\ref{thm: Polish space cts open surj image of cN} is spelled out in some detail.
\medskip

\noindent\textbf{AI Disclosure.} AI was not used during the research or writing of this work.

\section{Preliminaries}\label{sec: preliminaries}

In this section, we provide the necessary topological background for what follows. Rather than overwhelming the reader with a long list of standard topological definitions, we explain only those terms that are directly relevant to our discussion, referring to the relevant sections of \cite{Munkres} for the remaining definitions and results, and providing examples and exposition where appropriate.\footnote{As Paul Halmos observed, mathematics often struggles in ``matching precise definitions to
our intuitive preconceptions'' \cite[p.63]{Halmos1985-automathography}. Standard topological notions such as \textit{continuity} and \textit{connectedness} often elude genuine understanding until one confronts examples that defy na\"ive intuition---such as a function that is discontinuous at every rational point and continuous at every irrational point, or a curve composed of two disjoint parts that is nevertheless considered connected. In this regard, we refer the interested reader to \cite{Steen-Seebach-1978} for a rich collection of such pathological examples.}

Topology can be thought of as the study of \textit{continuous functions} between \textit{topological spaces}, which are pairs $(X,\tau)$, where $X$ is a set and $\tau$ is a collection of subsets of $X$ called a {\it{topology}} on $X$. The sets in $\tau$, called the {\it{open sets}}, satisfy the conditions that: $\es\in \tau$, $X \in \tau$, any union of sets in $\tau$ is in $\tau$, and any finite intersection of sets in $\tau$ is in $\tau$.

Given topological spaces $(X,\tau)$ and $(Y,\tau')$, a function $f:X\to Y$ is said to be {\it{continuous}} if, for every $S \in \tau'$, $f^{-1}(S) \in \tau$. Furthermore, a function $g: X\to Y$ is said to be a \textit{homeomorphism} if it is a bijection, and both $g$ and $g^{-1}$ are continuous \cite[\S18]{Munkres}. Topological spaces $X$ and $Y$ are said to be \textit{homeomorphic}, denoted $X\cong Y$, if there exists a homeomorphism between them. We can informally regard ``homeomorphic'' as meaning ``topologically identical''.

We use $\bR$ to denote the set of (ordinary) \textit{real numbers}, equipped with the \textit{standard topology}. The open sets can be generated by taking (possibly infinite) unions of open intervals of the form $(a,b) = \set{x\in \bR\mid a<x<b}$, where $a,b \in \bR$.\footnote{Here, $\bR$ can be written as the union of intervals of the form $(-a,a)$ for $a\in \bR$, and $\es$ can be thought of as the special case of $(a,b)$ where $a=b$.} With this topology, $\bR$ is \textit{homeomorphic}  to any nonempty open interval, considered as a subset of $\bR$ equipped with the \textit{subspace topology} \cite[\S16]{Munkres}. So, for example, $\bR$ is homeomorphic to $(0,1), (1,\infty),$ and $(0,\infty)$---a fact we shall take advantage of in Section \ref{sec: Reducing R to Baire space}.

$\bR$, as a topological space, has several favorable properties. For instance, $\bR$ satisfies many properties that capture the intuitive notion of being ``connected'': $\bR$ is \textit{connected}, \textit{locally connected}, and every point in $\bR$ is a \textit{strong cut point}, i.e., removing it leaves precisely two \textit{connected components} \cite[\S25]{Munkres}.\footnote{Let $X$ be a topological space. $X$ is \textit{connected} if it cannot be written as a disjoint union of nonempty open subsets; and a subspace $Y\subseteq X$ is \textit{connected} in $X$ if it cannot be written as a union $Y=A\cup B$ of disjoint nonempty sets (in the subspace topology), neither of which contains a limit point of the other \cite[\S23]{Munkres}. $X$ is \textit{locally connected} if for any $x\in X$ and any open set $U\ni x$, there is a connected open set $C\ni x$ such that $C\subseteq U$ \cite[\S25]{Munkres}.}

Furthermore, $\bR$ is \textit{metrizable}, which roughly means that it is possible to define a notion of distance (a \textit{metric}) between points in $\bR$ that is compatible with its topological structure, in the sense that the metric induces the topology on $\bR$; see \cite[\S20]{Munkres} for precise definitions.\footnote{In fact, $\bR$ is \textit{completely metrizable}. This means that it is possible to define a metric $d(x,y)$ on $\bR$ that is \textit{Cauchy-complete}, i.e., every \textit{Cauchy sequence} in $\bR$ (sequence of points that are getting arbitrarily close, in some precise sense, with respect to the metric $d$) \textit{converges} to a point in $\bR$. For example, the standard metric on $\bR$ is complete \cite[Theorem 43.2]{Munkres}. This result plays a key role in the field of introductory real analysis. See \cite[\S43]{Munkres} for a treatment of complete metric spaces, and \cite{Abbott2015} for an accessible reference in introductory real analysis.\label{footnote: completely metrizable}} For example, the \textit{standard (Euclidean) metric} on $\bR$ is given by $d(x,y)=\abs{x-y}$; and, the standard topology on $\bR$ can be generated by taking \textit{open balls} $B_\varepsilon(a)=\set{x\in \bR \mid d(a,x)<\varepsilon}$ of radius $\varepsilon$ centered at $a$, for $\varepsilon \in\bR_{>0}$ and $a\in \bR$.

Finally, $\bR$ is \textit{separable}, meaning that $\bR$ contains a countably infinite \textit{dense} subset, i.e., a countably infinite subset whose \textit{topological closure} is $\bR$ \cite[\S17]{Munkres}.\footnote{Since the \textit{order-topology} on $\bR$ is the same as its standard topology \cite[\S14]{Munkres}, $A\subseteq\bR$ is dense in the topological sense, if and only if $A$ is dense in the standard ordering on $\bR$, i.e., for every $x,y\in \bR$, with $x<y$, there exists $a\in A$ such that $x<a<y$. $\bQ$, the set of rational numbers, is an example of a countably infinite dense subset of $\bR$.} 

These properties---connectedness, local connectedness, separability, metrizability, and every point a strong cut point---are collectively sufficient to characterize $\bR$ as a topological space:\footnote{There are other characterizations of $\bR$: As an ordered field, $\bR$ is characterized as the \textit{unique up to isomorphism (Dedekind-)complete ordered field}---where, \textit{Dedekind-complete} means that every nonempty subset of $\bR$ that is bounded above has a least upper bound. ($\bQ$, for example, fails to satisfy this, since any sequence of rationals converging to an irrational from below does not have a least upper bound in $\bQ$.) 
One way to think about Dedekind-completeness is that $\bR$ has no “holes,” while $\bQ$ does. This characterization is central in other areas of mathematics, such as real analysis (see, for example, \cite[\S1.3]{Abbott2015}). A characterization of the ordering of $\bR$ can be found in \cite[Theorem 2.30]{Rosenstein-linear-orderings-1982}.\label{footnote: R complete ordered field}}

\begin{theorem}[Direct consequence of \cite{Ward-1936}]\label{thm: top char of R}
$\bR$ is the unique (up to homeomorphism) connected, locally connected separable metric space, such that every point is a strong cut point.
\end{theorem}

The properties are not only collectively sufficient, but also individually necessary, in that Theorem \ref{thm: top char of R} fails if one of them is omitted (see, for example, \cite{Steen-Seebach-1978}). Because the properties listed in Theorem~\ref{thm: top char of R}, as well as ``completeness'' properties such as complete metrizability (see footnotes~\ref{footnote: completely metrizable} and \ref{footnote: R complete ordered field}), are satisfied by spaces homeomorphic to $\bR$, set theorists often refer to any such space, including $\bR$ itself, as \textit{the continuum}.

Despite this, much of modern descriptive set theory centers not on $\bR$ itself, but on the Baire space $\cN = {}^\omega\omega$, whose elements are often called ``reals,'' as noted in footnote~\ref{footnote: examples where Baire space is called reals}. To understand this substitution, we now define $\cN$ and examine its topological properties in comparison with $\bR$. When making the distinction is necessary, we will refer to elements of $\bR$ as \textit{ordinary reals} and elements of $\cN$ as \textit{set-theorist's reals}.

Given a (possibly infinite) collection $\set{X_i}_{i\in I}$ of topological spaces, the \textit{product topology} on the cartesian product $\prod_{i\in I} X_i$ can be defined as follows: $U\subseteq \prod_{i\in I} X_i$ is a basic open set in the product topology if $U=\prod_{i\in I} U_i$, where each $U_i$ is open in $X_i$ and all but finitely many $U_i$ are equal to $X_i$.\footnote{It is worth noting that if $I$ is infinite (like the case we are concerned with, where $I=\omega$) the product topology can differ from the \textit{box topology} which allows infinitely many $U_i$'s to be proper subsets of $X_i$. For the categorically inclined reader, we note that since the product topology on $\prod_{i\in I} X_i$ is the \textit{coarsest} (smallest by inclusion) topology for which all the projections $\prod_{i\in I} X_i\to X_j$, for $j\in I$, are continuous, $\prod_{i\in I} X_i$ equipped with the product topology is the \textit{categorical product} (see \cite[\S III.4]{MacLane-Categories}) in the category of topological spaces. For details on the product topology, and exposition on why it is preferred over the box topology, see \cite[\S19]{Munkres}.\label{footnote: product topology}} The open sets are arbitrary unions of these basic open sets.

\begin{definition}\label{def: Baire space} Let $\omega = \set{0,1,2,\ldots}$ denote the set of natural numbers. The \textit{Baire space}, $\cN={}^\omega\omega$, has underlying set ${}^\omega\omega$, the set of countable sequences of natural numbers (i.e., functions $\omega\to \omega$), or equivalently, the countable product $\prod_{n\in\omega}\omega$. Its topology is given by the product topology, where each constituent $\omega$ is given the \textit{discrete topology} (every subset of $\omega$ is open).
\end{definition}

Henceforth, elements of $\cN$ are listed as $a_0, a_1,a_2,\ldots$, where each $a_i$ is a sequence coded as a function, which is listed as $a_i(0),a_i(1),a_i(2),\ldots$, to avoid confusion or the use of two levels of subscripts.

In a similar manner to Theorem~\ref{thm: top char of R}, $\cN$ has its own topological characterization, which, as we will indicate in Section~\ref{sec: Reducing R to Baire space}, can be used to prove that some topological space is homeomorphic to $\cN$:

\begin{theorem}[Alexandrov--Urysohn; see Theorem 7.7 in \cite{Kechris-classical-DST-1995}]\label{thm: top char of N}
$\cN$ is the unique (up to homeomorphism) nonempty zero-dimensional Polish space for which all compact subsets have empty interior. (See \cite{Munkres,Kechris-classical-DST-1995} for relevant definitions.)
\end{theorem}

It is only natural to ask if the aforementioned favorable properties of $\bR$ carry over to $\cN$. $\cN$ is indeed separable and completely metrizable. But as far as connectedness is concerned, $\bR$ and $\cN$ are as different as can be: $\cN$ is \textit{totally disconnected}, meaning that its only connected components are singletons (one-element sets). As Andreas Blass wryly put it: ``What's the continuum? Well it's this Baire discontinuum.''\footnote{Remark made by Andreas Blass in conversation with the first author.} In Section~\ref{sec: favorable properties of Baire space}, we will see how these connectedness properties render $\bR$ relatively inflexible for the purposes of DST, while their absence makes $\cN$ more favorable to work with.

How do we know that $\cN$ is separable and completely metrizable?

\begin{enumerate}[topsep=0pt,itemsep=0pt,partopsep=1ex,parsep=1ex]
\item[(i)] To show separability, we let $^{<\omega}\omega$ denote the set of finite strings of natural numbers. Sets of the form $U_\sigma = \set{a\in \cN \mid a\succ \sigma}$, for $\sigma \in\ ^{<\omega}\omega$ (i.e., sequences in $\cN$ extending the finite string $\sigma$), form a base for the topology on $\cN$. Then, it is straightforward to see that the set of all eventually zero sequences of natural numbers, i.e., the set $\set{a\in \cN \mid \exists j\in \omega\ (\forall k\geq j\ (a(k)=0))}$, is a countably infinite dense subset of $\cN$, so $\cN$ is separable.

\item[(ii)] The map $d:\cN\times\cN \to \bR$ defined as $d(a,b) = \frac{1}{2^{n+1}}$, where $n$ is the least element of $\omega$ such that $a(n)\neq b(n)$ (with $d(a,a)=0$) is a complete metric on $\cN$.
 \end{enumerate}\smallskip

 There is a standard ordering on $\cN$, the \textit{lexicographical ordering}:\ for $a,b\in \cN$, we say $a<_{\text{lex}} b$ if either $a(0)<b(0)$ or \big($a(0)=b(0) \text{ and } a(1)<b(1)$\big) or \big($a(0)=b(0) \text{ and } a(1)=b(1) \text{ and } a(2)<b(2)$\big) or $\cdots$. The same set of all eventually zero sequences of natural numbers, i.e., $\set{a\in \cN \mid \exists j\in \omega\ (\forall k\geq j\ (a(k)=0))}$, is a countably infinite dense (in the sense of order) subset of $\cN$ in the lexicographical ordering.\footnote{As we will see in Section~\ref{sec: contrasting maps f and g}, the $<_{\text{lex}}$-order topology differs from the product topology on $\cN$ and the resulting spaces are not homeomorphic.}
 
In fact, these two properties, separability and complete metrizability---shared not only by $\bR$ and $\cN$ but also by many other spaces of interest---are precisely the ones that play a central role in DST:

\begin{definition}\label{def: Polish space} A \textit{Polish space} is a separable, completely metrizable topological space.
\end{definition}

Examples of Polish spaces include $\bR$, $\bR^n$ with the product topology, $\cN$, and \textit{Cantor space} $\mc{C}={}^\omega2$ (which, since $2=\set{0,1}$ is a subset of $\omega$, can be defined via the subspace topology $\mc{C}\subseteq\cN$).\footnote{In certain contexts, ``reals'' are identified with elements of $\mc{C}$. An example in the classical case is \cite[p.31]{Becker-1988}, where he mentions that the results still hold for $\cN$: ``In this paper we work with Cantor space $\mc{C}=2^\omega$; all of our results are also valid for $\omega^\omega$.'' For an example in the effective context (particularly, countable admissible ordinals), see \cite[p.9]{Arant-Kechris-Lutz-2024}: ``for any real $x\in 2^\bN$.'' In Section \ref{sec: favorable properties of Baire space}, we will define the projective hierarchy (Definition \ref{def: projective hierarchy}), where $\cN$ is a favorable choice over both $\bR$ and $\mc{C}$. In Section \ref{sec: reducing general case to Baire space}, we will see more results that hold for $\cN$ (Theorems~\ref{thm: Polish space cts open surj image of cN} and \ref{thm: perfect Polish space cts bij image of cN}) but not for $\mc{C}$.\label{footnote: Cantor space as reals}} More examples of Polish spaces can be found in \cite[\S3.A]{Kechris-classical-DST-1995}. 

These spaces also share another property that is of interest in this context: They have no \textit{isolated points}, i.e., there is no point $x$ such that the singleton $\set{x}$ is open. A topological space with no isolated points is said to be \textit{perfect}. More examples of perfect Polish spaces can be found in \cite[\S6.A]{Kechris-classical-DST-1995}.

Since Polish spaces are metrizable, they are \textit{Hausdorff} (also known as $T_2$), meaning that distinct points are separated by disjoint open sets. More precisely, if $X$ is a nonempty Polish space and $x,y\in X$ are distinct, then there exist disjoint open sets $U,V\subseteq X$ such that $x\in U$ and $y\in V$.\footnote{In fact, metrizable spaces (hence, Polish spaces) satisfy much stronger separation axioms; see, for example, \cite[\S31]{Munkres} and \cite[Theorem 32.2]{Munkres}.} This can be seen as follows: Given $d$, a metric on $X$, and distinct points $x,y\in X$, let $\varepsilon=\frac{d(x,y)}{2}$ and define $U=B_\varepsilon(x)$, i.e., the open ball of radius $\varepsilon$ centered at $x$, and similarly define $V=B_\varepsilon(y)$. Then, by definition, $x\in U, y\in V$, and $U\cap V=\es$. 

Hausdorff spaces have many nice properties; for instance, any sequence $(y_n)$ of points in a Hausdorff space $Y$ converges to at most one point \cite[Theorem 17.10]{Munkres}. However, there are several naturally occurring non-Hausdorff topological spaces throughout mathematics. 

For example, the \textit{cofinite topology} on an infinite set $S$ (where a subset $U$ of $S$ is open if and only if $S\setminus U$ is finite or $U=\es$) is $T_1$---meaning that singletons, hence finite sets, are closed---but not $T_2$. The cofinite topology occurs naturally in algebraic geometry as the \textit{Zariski topology} on the set of closed points of an \textit{algebraic variety}. In general, if $R$ is a ring (commutative with unit), then the \textit{Zariski topology} on $\Spec(R)$ (the set of prime ideals of $R$) is $T_0$---meaning that for every pair of distinct points, at least one of them is contained in an open set not containing the other---but not necessarily $T_1$. (A non-maximal prime ideal in $R$ corresponds to a point that is not closed.) We refer the interested reader to \cite{Eisenbud1995} for details on the Zariski topology and \cite{Hofmann-2003} for a brief discussion on naturally occurring non-Hausdorff topological spaces.

Thus, it is worth studying non-Hausdorff generalizations of Polish spaces, such as the well-behaved \textit{quasi-Polish spaces}, which we discuss in Section \ref{sec: reducing general case to Baire space}.

\section{Favoring the ``Baire discontinuum''}\label{sec: favorable properties of Baire space}

In this section, we will see how the connectedness of $\bR$ makes it inflexible in the context of DST, favoring the radical not-connectedness of $\cN$. So, the very properties that disqualify $\mathcal{N}$ as a substitute for $\mathbb{R}$ in other contexts turn out to be structural advantages for DST.

A defining feature of DST is that it measures the complexity of sets (subsets of a Polish space) in terms of the complexity of their definition. This is formalized by the \textit{projective hierarchy} (Definition~\ref{def: projective hierarchy}) on a given Polish space. We follow the notation in \cite[\S 2]{Martin-DST-Proj-Sets-1977}.

\begin{definition}\label{def: projection} For positive integer $k$ and $B\subseteq \cN^{k+1}$, the \textit{projection} of $B$ is the set \[\set{(a_1,\ldots,a_k) \in \cN^k \mid \exists b\in \cN\big( (a_1,\ldots,a_k, b)\in B \big)}.\]
\end{definition}

\begin{definition}\label{def: projective hierarchy}
We inductively define the {\textit{projective hierarchy}} on a Polish space $X$ (where $X$ is $\cN$ or $\cN^m$)\footnote{The projective hierarchy can be defined for more general Polish spaces by modifying the definition of projections to sets of the form $B\subseteq X\times Y$. (Our definition is for the special case of $X=\cN^k$ and $Y=\cN$.) In this case, a projection along $Y$ can be thought of as existential quantification over $Y$, as described in \cite[p.25]{Moschovakis-DST}. However, as we will see in Sections \ref{sec: Reducing R to Baire space} and \ref{sec: reducing general case to Baire space}, for the purposes relevant to us, $\cN^m$ suffices.} as follows. For all $n\in \omega$ and $A\subseteq X$:
\begin{align*}
    A\in \bSigma^1_0 &\iff A \text{ is open} \\
    A\in \bPi^1_n &\iff X\setminus A \in \bSigma^1_n \\
    A\in \bSigma^1_{n+1} &\iff \text{there is } B\in \bPi^1_n \text{ such that } A \text{ is the \textit{projection} of } B
    \end{align*}
\end{definition}

For orientation, it's useful to consider the logical complexity of the definition of a set in this hierarchy. First note that projection, according to Definition \ref{def: projection}, incorporates a reduction of dimension by one variable. The projection of $B \subseteq \cN^{k+1}$ will be a set $A \subseteq \cN^k$. This isn't just a notational convenience: the reduction in dimension arises from the logical structure of the definition of projection. The core logical operation in the definition of projection above is the existential quantifier, which is essentially a reduction in dimension.\footnote{The interaction between projection and dimension can be further clarified by noting that in some mathematical contexts, projections can appear in the form of existential quantification, wherein the $\exists$ quantifier turns an $n+1$ place relation into an $n$ place relation. For example, let $R(x,y)\subseteq\omega\times\omega$ be given by $x=2y+1$, and consider the projection $S(x) = ``\exists yR(x,y)$'', i.e., $S(x) = \set{x\in \omega \mid \exists y\in \omega \big(R(x,y) \big)}$. Then $S(x)$ holds if and only if $x$ is an odd number.}

The complexity of the sets in the projective hierarchy is straightforwardly displayed in the \textit{quantifier complexity} (formally defined using the \textit{analytical hierarchy}) with additional function parameters. We've noted that projection corresponds logically to existential quantification, and of course complementation corresponds to logical negation ($\neg$). Since the universal quantifier $\forall x$ can be interchanged with the operator $\neg \exists x\neg$, a sequence of projections and complementations will be marked by strings of alternating existential and universal quantifiers.\footnote{For instance, suppose we have a set defined by a formula $\Phi$---with a sufficient number of free (function) variables, as the $a_i$'s are officially functions. Since, for example, $\neg\exists x_1\neg\exists x_2 \Phi$ is equivalent to $\forall x_1 \exists x_2\Phi$, if $\Phi$ defines an open set, the formula $\forall x_1\exists x_2\Phi$ would define a $\bPi^1_2$ set. Similarly, the $\bSigma^1_n$ (respectively, $\bPi^1_n$) sets are those that can be defined by a string of $n$ alternating function quantifiers beginning with an existential (respectively, universal) quantifier $\exists a_1\forall a_2 \ldots Q_na_n \Phi$, with $\Phi$ defining an open set, and $Q_n$ either $\forall$ or $\exists$. For a more precise and detailed version of these brief orienting remarks about quantifier complexity and the definition of the projective hierarchy, see \cite[Ch.IV]{Hinman-RecursionHierarchies-1978}, and see earlier chapters for definitions.}

As we unravel the definition of the projective hierarchy---saying that $\bSigma^1_n$ sets are projections of $\bPi^1_{n-1}$ sets of higher dimension, which are in turn complements of $\bSigma^1_{n-1}$ sets, which again are projections of $\bPi^1_{n-2}$ sets of a yet higher dimension, and so on, eventually down to open sets---we find ourselves needing to go to higher and higher dimensions.

This is a feature of the projective hierarchy that $\cN$ is well-equipped to handle, since $\cN$ is homeomorphic to $\cN^2$, i.e., its product with itself, whereas $\bR$ is not homeomorphic to $\bR^2$. For instance, the map $\cN \to \cN^2$ given by $(a(i)) \mapsto (a(2i), a(2i+1))$ is a homeomorphism. Given two (and similarly $k$) variables ranging over $\cN$, we can treat them as a single variable by interleaving the sequences. This argument generalizes to show that $\cN$ is homeomorphic to $\cN^k$ for any positive integer $k$. Martin describes this map as preserving ``everything of importance in classical and effective descriptive set theory'' \cite[p.786]{Martin-DST-Proj-Sets-1977}.  Thanks to this feature of  $\cN$, ``the notion of dimension plays no role in descriptive set theory. Thus, although our theory deals with finite
product spaces $X^n$, we shall often state definitions and theorems simply for the case $\cN$, and no generality is lost in doing so'' \cite[p.786]{Martin-DST-Proj-Sets-1977}. 
In contrast, we must modify the definition of $\bSigma^1_n$, by allowing more complicated sets in $\bSigma^1_0$, for the ordinary reals (or even Cantor space, $\mc{C}$, defined in Section~\ref{sec: preliminaries}) \cite[Remark 2.12]{Martin-DST-Proj-Sets-1977}. Martin also explicitly contrasts $\cN$ with $\bR$, in this regard, noting it as a reason to regard the ordinary reals as ``slightly awkward to use'' \cite[p.785]{Martin-DST-Proj-Sets-1977}: ``the fact that the simplest bijection of the [ordinary] reals onto the plane is slightly more complicated than a homeomorphism is one reason we do not base [DST] on the [ordinary] reals.'' \cite[p.786]{Martin-DST-Proj-Sets-1977}.


\section{Another topological difference}\label{sec: another topological difference}

Before we consider in more detail the ways in which $\cN$ and $\bR$ are topologically similar in ways relevant to descriptive set theorists, it's worth emphasizing that there are many non-trivial topological differences between them that can be profoundly significant in other areas of mathematics. We have already seen this in the form of the many ``connectedness'' properties of $\bR$ (as in Theorem~\ref{thm: top char of R}) and the total disconnectedness of $\cN$ prompting Blass's ``Baire discontinuum'' remark (see Section~\ref{sec: preliminaries}). We now compare these spaces (and Cantor space $\mc{C}$, defined in Section~\ref{sec: preliminaries}) in terms of ``compactness'' properties. 

First note that $\mc{C}$ is itself \textit{compact} \cite[\S26]{Munkres}, $\bR$---while not compact---is \textit{locally compact} \cite[\S29]{Munkres}, and $\cN$ is not even locally compact. On a related note, one may ask how many compact sets are needed to cover each of these spaces. Since $\mc{C}$ is itself compact, it needs exactly one. $\bR$ can be covered by (i.e., written as a union of) countably many compact subsets---for example, $\bR=\bigcup_{n\in\bZ_{>0}} [-n,n]$---hence, it is $\sigma$\textit{-compact}. Even this is not true of $\cN$, i.e., $\cN$ is not a countable union of compact subsets, hence uncountably many compact subsets are needed to cover $\cN$. To describe precisely how many compact sets are needed, we recall the dominating number.

For $f,g \in \cN$, define $f\leq^*g$ to mean $f(n)\leq g(n)$ for all sufficiently large $n$. A family $D\subseteq\cN$ is \textit{dominating} if every $f\in\cN$ satisfies $f\leq^*g$ for some $g\in D$. The \textit{dominating number} $\mathfrak d$ is the least cardinality of a dominating family. No countable family $\set{g_0,g_1,\ldots}$ is dominating: the function $h(n)=1+\max_{i\leq n}g_i(n)$ is not eventually bounded by any $g_i$. Thus, $\aleph_1\leq\mathfrak d\leq 2^{\aleph_0}$, making $\mathfrak d$ one of the \textit{cardinal characteristics of the continuum}. See \cite{Blass-cardinal-char-2010} for more details on such cardinals.

Every nonempty compact set $K\subseteq\cN$ is pointwise bounded: each coordinate projection of $K$ is a compact subset of the discrete space $\omega$, hence finite, so $b_K(n)=\max\{a(n):a\in K\}$ is well-defined. If a collection of compact sets covers $\cN$, then their pointwise bounds $b_K$ form a dominating family, so at least $\mathfrak d$ compact sets are needed to cover $\cN$. In fact, $\mathfrak d$ compact sets suffice \cite[Exercise 4.11]{Kechris-classical-DST-1995}. In particular, since $\mathfrak d>\aleph_0$, $\cN$ is not $\sigma$-compact.

Since $\cN$ is homeomorphic to the space of irrationals $\mc{I}=\bR\setminus\bQ$ (see Section~\ref{sec: irrationals map f}), precisely $\mathfrak d$ compact subsets are also needed to cover $\mc{I}$.\footnote{We refer the interested reader to \cite{Hechler-1974,Hechler-1975} for other natural occurrences of $\mathfrak d$.}


This $\sigma$-compactness difference between $\cN$ and $\bR$ puts them on different sides of a deep dichotomy. Indeed, it is a classical result of Hurewicz that a Polish space contains a closed subset that is homeomorphic to $\cN$ if and only if it is not $\sigma$-compact (see \cite[Theorem 7.10]{Kechris-classical-DST-1995}). In light of this result, Anush Tserunyan was moved to describe $\cN$ as the ``canonical obstruction to being $\sigma$-compact'' \cite[p.8]{Tserunyan-notes}.


\section{Reducing the case of $\mathbb{R}$ to $\mathcal{N}$}\label{sec: Reducing R to Baire space}

So $\bR$ is similar in certain ways to $\cN$, and $\cN$ is more fruitful for the descriptive set theorist than the ``awkward'' (for DST) $\bR$. But we've also seen that the two structures are different in mathematically significant ways. Indeed, since $\bR$ is connected, and $\cN$ isn't, they can't be homeomorphic. So how can we justify the common practice in DST to refer to elements of $\cN$ as ``reals''? As this section and Section \ref{sec: reducing general case to Baire space} will explain, they are ``sufficiently alike'' that---in some precise sense that will be explained later in the section---any DST fact proven about $\cN$ carries over to $\bR$ (and many other spaces too).\footnote{This discussion draws from Blass's comments in \cite[p.397]{Blass-cardinal-char-2010}, as well as remarks made by Blass in conversation with the first author. More general cases beyond $\bR$ will be addressed in Section \ref{sec: reducing general case to Baire space}.}

A key observation is that there are many ``nice'' bijections between $\cN$ and ``almost all'' of $\bR$, or other spaces that may be more obviously thought of as ``the reals,'' such as $(1,\infty), (0,1)$ and $[0,1)$. Somewhat more precisely, we present two ways in which an appropriate countable subset of $X$ (where $X$ is one of the aforementioned spaces) can be ignored to produce a continuous bijection from $\cN$.\footnote{Indeed, countable subsets are well worth ignoring for the purposes of DST. For instance, given an uncountable Polish space $X$ and a \textit{Borel set} $B\subseteq X$, if $B$ is modified to $B'$ by adding or removing countably many points, then $B'$ is still Borel (in fact, of the same level as $B$ in the Borel hierarchy, except for very low levels). Moreover, if $A$ is a countable subset of $X$, then $A$ is \textit{meager} (one way of being topologically insignificant) and \textit{measure zero} (measure-theoretically insignificant). While the measure zero result is briefly described in footnote \ref{footnote: Lebesgue Measure}, we refer the interested reader to \cite{Oxtoby-1980, Kechris-classical-DST-1995, Stein-Shakarchi-3} for details---the first of which describes connections between (Lebesgue) measure and (Baire) category.} The first is occasionally mentioned---though rarely explicitly---in the DST literature; the second appears even less frequently, as we explain later when contrasting these maps.

\subsection{Irrationals}\label{sec: irrationals map f}

The first construction is given by the \textit{continued fraction expansion}. Any point $x\in (0,\infty)$ can be expanded as a continued fraction:
\begin{equation*} x = n_0(x) + \frac{1}{n_1(x) + \frac{1}{n_2(x) + \frac{1}{\tiny\cdots}}}, \end{equation*}
where each $n_i(x)$ is in $\omega$ \cite[\S X]{Hardy-Wright}.\footnote{For the purposes of number theory, this discussion can be modified appropriately for $x\in \bR$ by allowing some of the $n_i$'s to be negative integers. For our purposes, it is enough that $\bR$ is homeomorphic to each of $(0,1), (1,\infty),$ and $(0,\infty)$, as discussed in Section \ref{sec: preliminaries}.} It is a standard number-theoretic result that every irrational number can be expressed uniquely as an infinite continued fraction. 

Let $\mc{I} = \bR\setminus\bQ$ denote the set of irrational numbers. Then, for any $x\in (0,\infty)\cap \mc{I}$, the continued fraction algorithm produces a unique sequence in $\cN$, given by $a(i)=n_i(x)$ \cite[Theorem~170]{Hardy-Wright}. This defines a function $\hat{f}:(0,\infty)\cap\mc{I} \to \cN$. Furthermore, if for $x,y \in (0,\infty)\cap \mc{I}$, we have $\hat{f}(x)=\hat{f}(y)$, then since two infinite continued fractions with the same \textit{partial quotients} (i.e., the same $n_i$'s) are equal \cite[Theorem~169]{Hardy-Wright}, $x=y$, whence $\hat{f}$ is injective. Note, however, that $\hat{f}$ is not surjective---since, for example, the constant zero sequence in $\cN$ is not mapped to. We modify $\hat{f}$ in the following manner to account for this issue.

We restrict the domain and codomain of $\hat{f}$ to produce \[f: (1,\infty)\cap\mc{I} \to {}^\omega(\omega\setminus\set{0}),\] which we first argue is a bijection. Since $x>1$ and irrational, $f(x)$ is still well-defined and now in ${}^\omega(\omega\setminus\set{0})$. Furthermore, since $f$ is a restriction of $\hat{f}$, $f$ is injective. Finally, given $a\in {}^\omega(\omega\setminus\set{0})$, let $x=a(0)+\frac{1}{a(1)+\cdots}$, which is in $(1,\infty)\cap\mc{I}$, since every infinite continued fraction converges to an irrational \cite[Theorem~166]{Hardy-Wright}. It is not hard to show that $f(x)=a$, so $f$ is surjective \cite[\S X]{Hardy-Wright}. Moreover, when $(1,\infty)\cap\mc{I}$ is considered as a subspace of $\bR$ and ${}^\omega(\omega\setminus\set{0})$ as a subspace of $\cN$, $f: (1,\infty)\cap\mc{I} \to {}^\omega(\omega\setminus\set{0})$ is a homeomorphism.\footnote{See, for example, \cite[Proposition 7.5]{Bertsekas-Shreve-StochasticOptimalControl-1978} for details on the continued fraction map, and \cite[Theorem 1.1]{Miller-DST-and-Forcing-2016} for an alternative proof that $\cN$ is homeomorphic to $\mc{I}$.}

Moreover, the map $t:\cN\to {}^\omega(\omega\setminus\set{0})$ given by $(a(i))_{i\in\omega} \mapsto (a(i)+1)_{i\in\omega}$ is a homeomorphism, so the composite $t^{-1}\circ f : (1,\infty)\cap\mc{I} \to \cN$ is a homeomorphism. Hence, up to a countable set, namely $(1,\infty)\cap\bQ$, $\cN$ and $(1,\infty)$ are topologically the same.

By composing with suitable homeomorphisms, $(1,\infty)$ can be replaced with $(0,1),(0,\infty), \bR$, etc. For example, composing $f^{-1}\circ t$ with the map $x\mapsto 1/x$ produces a homeomorphism $\cN \to (0,1)\cap \mc{I}$. Hence, up to a countable set, namely $(0,1)\cap\bQ$, $\cN$ and $(0,1)$ are also topologically the same.

\subsection{Non-dyadic rationals}\label{sec: non-dyadics map g}

The second construction is a continuous bijection \[g: \cN \to [0,1).\] We map $a \in \cN$ to the number whose binary expansion is a binary point, $a(0)$ ones, a zero, $a(1)$ ones, a zero, and so on, written in binary representation as:
\[ a\mapsto 0.\underbrace{1\ldots1}_{a(0)\text{-many}}0\underbrace{1\ldots1}_{a(1)\text{-many}}0\underbrace{1\ldots1}_{a(2)\text{-many}}0\cdots \]
It is not hard to show that $g$ is bijective; one can, for example, explicitly construct its inverse and easily verify that this construction is correct. It is even easier to see that $g$ is not a homeomorphism: $\cN$ is totally disconnected and $[0,1)$ is an interval, hence connected. Notably, $g^{-1}$ is continuous everywhere except at the dyadic rationals (rationals whose denominators are a power of $2$) between 0 and 1, denoted $\bZ[\frac{1}{2}]\cap(0,1)$, which again is a countable set. Hence, $g:\cN\to [0,1)$ is a continuous bijection whose inverse is discontinuous only on a countable set---which suffices for the purposes of descriptive set theory.

For completeness, we note that $g$ can be modified to a genuine homeomorphism, in the following manner. Why we may choose not to do so will be made more clear in the following section. Let $\mc{M} = \cN\setminus\set{a\in \cN \mid \exists j\in \omega\ (\forall k\geq j\ (a(k)=0))}$ denote the subspace of all non-eventually zero sequences in $\cN$. Then restricting the domain and codomain of $g$ produces a homeomorphism $\tilde{g}: \mc{M} \to [0,1)\setminus\bZ[\frac{1}{2}]$. It now suffices to show that $\cN\cong \mc{M}$. This can be done by showing that $\mc{M}$ satisfies the topological characterization of $\cN$ in Theorem~\ref{thm: top char of N}. Using Section \ref{sec: irrationals map f}, one could instead show that $(0,1)\setminus \bQ \cong (0,1)\setminus\bZ[\frac{1}{2}]$.\footnote{Since $(0,1)\cap\bQ$ and $(0,1)\cap\bZ[\frac{1}{2}]$ are countable dense linear orders with no end points, they are order-isomorphic (see \cite[Theorem 2.8]{Rosenstein-linear-orderings-1982}). This map extends uniquely to an order-isomorphism $(0,1)\to (0,1)$ sending $(0,1)\cap\bQ$ to $(0,1)\cap\bZ[\frac{1}{2}]$ (see \cite[\S2.4]{Rosenstein-linear-orderings-1982}). Since the order-topology on $\bR$ is the same as its standard topology \cite[\S 14]{Munkres}, this map is a homeomorphism. Hence, the restriction $(0,1)\setminus \bQ\to (0,1)\setminus\bZ[\frac{1}{2}]$ is a homeomorphism.\label{footnote: extending countable DLOs}}

\subsection{Contrasting $f$ and $g$}\label{sec: contrasting maps f and g}

As mentioned earlier, neither of these constructions is typically mentioned in the classical DST literature, though $f$ seems to be more common.\footnote{The continued fraction map is discussed in detail in \cite[Proposition 7.5]{Bertsekas-Shreve-StochasticOptimalControl-1978}, without proof in \cite{Jech-set-theory-1997}, as Exercise 3.4(ii) in \cite{Kechris-classical-DST-1995}, and as the following remark in \cite{Moschovakis-DST}: ``One can show that $\cN$ is homeomorphic with the set of irrational numbers topologized as a subspace of $\bR$. The proof appeals to some basic properties of continued fractions and does not concern us here---it can be found in any good book in number theory, for example \cite{Hardy-Wright}. Although we never use this result, we will find it convenient to call members of $\cN$ irrational.'' \cite[p.12]{Moschovakis-DST}. The latter construction appears in Blass's comments on the meaning of ``continuum'' in \cite[p.397]{Blass-cardinal-char-2010}. The first and a similar construction to the second appear as Exercises 1.13 and 1.12 in \cite{Marker-DST-notes}. It is also worth noting that a different homeomorphism $\cN\cong \mc{I}$ is discussed in \cite[Theorem 1.1]{Miller-DST-and-Forcing-2016} and $\cM\cong((0,1)\setminus\bZ[\frac{1}{2}])$ is discussed in \cite{Deiser-Baire-Dyadic-2009}.\label{footnote: second identification more common}} 
This is possibly because DST usually just needs the existence of such an identification rather than any specific correspondence (as highlighted by Moschovakis's remark in footnote~\ref{footnote: second identification more common}). Indeed, to a descriptive set-theorist, what matters is a homeomorphism up to countably many points, so that the topology is essentially the same. We still believe it is worth contrasting these maps, and argue---alongside Andreas Blass\footnote{Remark made by Andreas Blass in conversation with the first author.}---that if the actual correspondence did matter, then $g$ would likely have been preferred.

$g: \cN \to [0,1)$ is monotonic with respect to the lexicographical ordering on $\cN$ (defined in Section \ref{sec: preliminaries}): $a<_{\text{lex}} b$ if and only if $g(a)<g(b)$. Indeed, if $a<_{\text{lex}} b$, then there is a least index $i$ where $a(i)<b(i)$, and for all $j<i$, $a(j)=b(j)$. This means that the binary expansions of $g(a)$ and $g(b)$ agree through the first $i$ zeros and the following $a(i)$ ones. At the next digit, $g(a)$ has a $0$ while $g(b)$ has a $1$, whence $g(a)<g(b)$. The other direction follows by linearity of the order. So, if one wished to use the standard ordering on the ordinary reals, they may still reduce to $\cN$ and use $<_{\text{lex}}$ instead, under this identification.

It is worth noting that since the map $g : ({}^\omega\omega, <_{\text{lex}}) \to ([0,1),<)$ is an order isomorphism, it is a homeomorphism when ${}^\omega\omega$ has the $<_{\text{lex}}$-order topology. Composing with the map $x \mapsto \frac{x}{1-x}$ gives a homeomorphism onto the closed subspace $[0,\infty)$ of $\bR$. So, ${}^\omega\omega$ with the $<_{\text{lex}}$-order topology has dimension 1 and is connected. In particular, since $\cN$ is zero-dimensional and totally disconnected, this gives two ways to see that ${}^\omega\omega$ with the $<_{\text{lex}}$-order topology is not homeomorphic to $\cN$.

In contrast with $g$, $f$, as currently defined, is not monotonic with respect to $<_{\text{lex}}$. The general behavior of continued fractions is more complicated: If the continued fraction expansions of $x$ and $y$ first differ at index $i$, then $x<y$ if and only if $n_i(x)<n_i(y)$ when $i$ is even, and $n_i(x)>n_i(y)$ when $i$ is odd \cite[Theorem 152]{Hardy-Wright}. Hence, $f$ can be thought of as monotonic with respect to an order given by alternating the lexicographic and reverse lexicographic orders, which is harder to understand.

Moreover, $g$ behaves well with measure, insofar as descriptive set theory is concerned. A standard measure on the interval $[0,1)$ is the \textit{Lebesgue measure}, which, in many ways, is a natural measure.\footnote{For a formal and self-contained treatment of the Lebesgue measure on $\bR^n$, which is beyond the scope of this article, see \cite[Ch.1]{Stein-Shakarchi-3}. In the case of $\bR$, the intuition behind (Lebesgue) measurability comes from the ordinary notion of length: the measure of a point is zero and that of an open interval $(a,b)\subseteq \bR$ is its \textit{length}, denoted $m(a,b)=b-a$. Since any open subset $U$ of $\bR$ can be written uniquely as a disjoint union of open intervals \cite[Theorem 1.3]{Stein-Shakarchi-3}, the measure of $U$ can be defined as the sum of the lengths of these open intervals. Moreover, the Borel sets---defined as the smallest collection of subsets of $\bR$ that includes the open sets and is closed under taking complements and countable unions---are measurable \cite[Ch.1, \S3]{Stein-Shakarchi-3}. (For instance, since singletons are closed, hence Borel, all countable sets are Borel, hence measurable.) A subset $A$ of $\bR$ has \textit{measure zero} if it can be covered by arbitrarily small unions of open intervals: for any $\varepsilon>0$, there are open sets $U_i$ with $\sum_i m(U_i)<\varepsilon$ such that $A\subseteq \bigcup_i U_i$. (For instance, countably infinite sets have measure zero.) More generally, a set $A \subseteq \bR$ is \textit{measurable} if and only if it differs from a Borel set by a set of measure zero \cite[Corollary 3.5]{Stein-Shakarchi-3}.\label{footnote: Lebesgue Measure}} While there does not seem to be an equally natural way to define a measure on $\cN$, under the map $g$, the Lebesgue measure on $[0,1)$ corresponds to the product measure on $\cN$ obtained from the measure on $\omega$, where a point $n\in \omega$ is given measure $\frac{1}{2^{n+1}}$. (In some sense, this is the best one can do. Since we wish to have total measure 1, there is no way to assign a uniform measure on $\omega$, as even the smallest uniform assignment will have infinite total measure.) Again, since the set $\bZ[\frac{1}{2}]\cap(0,1)$ that $g^{-1}$ is discontinuous on is countable, it has Lebesgue measure 0, i.e. it is measure-theoretically insignificant.

In contrast, while there is a description of a measure on $(0,1)$ that is preserved by $f$, it is more complicated to understand than the Lebesgue measure. In particular, $f$ preserves the \textit{Gauss measure} $\mu_G$ on $(0,1)$ given by
\[ \mu_G(A) = \frac{1}{\log2}\int_A \frac{1}{1+x}dx, \]
for any Borel set $A\subseteq (0,1)$ \cite[Lemma 3.5]{Einsiedler-Ward-Ergodic-Theory-2011}.

Since a homeomorphism up to countably many points is sufficient for the purposes of Baire Category, both $f$ and $g$ preserve Baire Category properties of subsets (see \cite[Ch.8]{Kechris-classical-DST-1995} or \cite{Oxtoby-1980}). On the other hand, since $f$ is a homeomorphism, $f$ preserves Borel sets at all levels of the Borel hierarchy \cite[Ch.11B]{Kechris-classical-DST-1995}, while $g$, whose inverse has countably many discontinuities, preserves Borel sets at all levels beyond the second.

One might wonder why we choose not to modify $g$ to the homeomorphism $\tilde{g}:\mc{M}\to[0,1)\setminus\bZ[\frac{1}{2}]$, as defined in Section~\ref{sec: non-dyadics map g}, and compose it with a suitable homeomorphism $\mc{N}\to \mc{M}$ (see, for example, footnote~\ref{footnote: extending countable DLOs}). While one could argue, using the rest of this section, that for the purposes of DST, we only require a homeomorphism up to countably many points, this modification may seem like an improvement if it doesn't take away from any of the properties of $g$. Unfortunately, we lose monotonicity with respect to $<_{\text{lex}}$.

Since $\cN$ has a least element (namely, the constant zero sequence), whereas $\cM$ has neither a greatest nor a least element, there is no order-preserving or order-reversing bijection from $\cN$ to $\cM$, hence no monotonic (with respect to $<_{\text{lex}}$) homeomorphism from $\cN$ to $\cM$.\footnote{Of course, stringing along the maps in footnote~\ref{footnote: extending countable DLOs} and $f$, one could produce a map $\mc{N}\to \mc{M}$ under which $<_{\text{lex}}$ on $\mc{M}$ corresponds to the aforementioned alternating lexicographic and reverse lexicographic order on $\cN$. However, as also mentioned earlier, this is more complicated to understand.}

\section{Reducing the general case to $\cN$}\label{sec: reducing general case to Baire space}

Although the reals are the main interest of DST, the subject applies more generally to all Polish spaces (see Definition~\ref{def: Polish space}) sharing the relevant topological properties of $\bR$. Recent fruitful investigations have also shown that the subject applies to quasi-Polish spaces: A \textit{quasi-Polish space} is a second-countable, completely \textit{quasi-metrizable} topological space, where a quasi-metric is a generalization of a metric that is not required to obey symmetry. Quasi-Polish spaces were introduced by de Brecht as well-behaved non-Hausdorff generalizations of Polish spaces that share many of their descriptive set-theoretic properties \cite{deBrecht-QuasiPolish-2013}.\footnote{A self-contained treatment of quasi-Polish spaces, approached from a more classical DST perspective, is provided in \cite{Chen-QuasiPolish-2018}.} In this section, we explore why, for many purposes, it still suffices to study $\cN$ even in these more general contexts. 

In the general topological setting, given spaces $X$ and $Y$, if $f: X\to Y$ is a continuous open (or closed) surjection---hence, a \textit{quotient map}---then the topology on $Y$ is necessarily the quotient topology induced by the topology on $X$ \cite[\S22]{Munkres}.\footnote{Given a topological space $X$, a set $A$, and a surjection $p:X\to A$, the \textit{quotient topology} on $A$ is the largest (by containment) topology on $A$ such that $p:X\to A$ is continuous. More explicitly, $U\subseteq A$ is open in the quotient topology if $p^{-1}(U)$ is open in $X$. For a detailed discussion, see \cite[\S22]{Munkres}.} This result roughly means that if one is interested in studying a class of spaces $Y$, which can all be written as open (or closed) quotients of a restricted class of spaces $X$, then it suffices to study the spaces $X$ together with the induced equivalence relations on them.\footnote{Given a surjection $f: X\to Y$, the induced equivalence relation $E_f$ on $X$ is defined, for every $x,x'\in X$, as $x\ E_f\ x'$ if and only if $f(x)=f(x')$.}

In the context of classical DST, where the relevant spaces $Y$ are Polish spaces, the situation is especially well-behaved, in the sense that there is a single space $X$---the familiar hero $\cN$ of this folktale---from which all nonempty Polish spaces can be obtained as continuous open images:

\begin{theorem}[7.14 in \cite{Kechris-classical-DST-1995}]\label{thm: Polish space cts open surj image of cN}
If $Y$ is a nonempty Polish space, then there is a continuous open (or closed \cite{Engelking-1969}) surjection $f : \cN \to Y$.\footnote{One can obtain such a continuous open surjection by following the steps in \cite[7.14]{Kechris-classical-DST-1995}. This argument is also spelled out in \cite{Marker-DST-notes}.}
\end{theorem}

This means that an arbitrary nonempty Polish space $Y$ is necessarily endowed with a quotient topology coming from $\cN$. Since Theorem~\ref{thm: Polish space cts open surj image of cN} holds for \textit{every} nonempty Polish space $Y$, it is an example of what is called a \textit{universality property}.\footnote{Another example of a universality property for $\cN$ is \cite[Theorem~7.9]{Kechris-classical-DST-1995}: For any Polish space $Y$, there is a closed subspace $F\subseteq \cN$ and a continuous bijection $h:F\to Y$ \cite[Theorem~7.9]{Kechris-classical-DST-1995}. This, along with the fact that any nonempty closed subset $F\subseteq \cN$ is a \textit{retract} of $\cN$---that is, there is a continuous surjection $g:\cN\to F$ such that $g(x)=x$ for any $x\in F$ \cite[Proposition~2.8]{Kechris-classical-DST-1995}---can be used to obtain, for any nonempty Polish space $Y$, a continuous surjection $f=h\circ g: \cN\to Y$. This result, as currently stated, is weaker than Theorem~\ref{thm: Polish space cts open surj image of cN}, as, a priori, there is no reason for this map to be open. We also note that the map $h:F\to Y$ is not generally a homeomorphism since closed subsets of $\cN$ are necessarily totally disconnected.}


 In footnote~\ref{footnote: Cantor space as reals}, we mentioned that in certain contexts, reals are identified with elements of Cantor space $\mc{C}$. Since $\mc{C}$ is compact, if there is a continuous surjection $f: \mc{C} \to Y$, then its image $f(\mc{C})=Y$ must also be compact \cite[Theorem~26.5]{Munkres}, hence Theorem~\ref{thm: Polish space cts open surj image of cN} does not hold if $\cN$ is replaced by $\mc{C}$. (For instance, since neither $\bR$ nor $\cN$ are compact, there is no continuous surjection from $\mc{C}$ onto either.) There is an appropriate universality property for $\mc{C}$: If $Y$ is a nonempty compact metrizable space, then there is a continuous surjection $f: \mc{C} \to Y$ \cite[Theorem~4.18]{Kechris-classical-DST-1995}.


If, instead of arbitrary Polish spaces, we are specifically concerned with perfect Polish spaces, then $\cN$ gives the following characterization:

\begin{theorem}[7.15 in \cite{Kechris-classical-DST-1995}]\label{thm: perfect Polish space cts bij image of cN} Suppose $Y$ is a nonempty Polish space. Then $Y$ is perfect if and only if there is a continuous bijection $f : \cN \to Y$.
\end{theorem}

By virtue of the two preceding theorems, in proving that a property holds of an arbitrary nonempty (perfect) Polish space, it suffices to show that said property is true of $\cN$ and preserved by continuous images.\footnote{A recent application of Theorem~\ref{thm: perfect Polish space cts bij image of cN} is in \cite[Theorem 3.1.1]{Frisch-Kechris-Shinko-Vidnyánszky-2023}: ``Since for every perfect Polish space $Y$, there is a continuous bijection $\bN^\bN\to Y$, we can assume $Y=\bN^\bN$.'' Since, by the Cantor--Bendixson theorem \cite[Theorem 6.4]{Kechris-classical-DST-1995}, any Polish space $X$ can be written uniquely as the disjoint union $X=P\cup C$ of a perfect set $P$ (the \textit{perfect kernel}) and a countable set $C$, Theorem~\ref{thm: perfect Polish space cts bij image of cN} can be applied to the perfect kernel of a nonempty Polish space. For an example, see B. Miller's proof of the Kechris--Solecki--Todorcevic $G_0$ dichotomy theorem \cite{Miller-Graph-Theory-DST-2012}.\label{footnote: use of cts bijection from N}}

$\cN$'s canonical character is robust enough to remain even in a natural generalization of Polish space. As mentioned in Section~\ref{sec: preliminaries}, there are several naturally occurring non-Hausdorff topological spaces throughout mathematics, so one may wish to study non-Hausdorff generalizations of Polish spaces, such as the aforementioned quasi-Polish spaces. Quasi-Polish spaces are not required to satisfy any separation axioms beyond $T_0$, hence distinct points need not be separated by disjoint open sets as in the Hausdorff setting of Polish spaces. What is remarkable is that quasi-Polish spaces not only satisfy the analog of Theorem~\ref{thm: Polish space cts open surj image of cN}, but are characterized by it:

\begin{theorem}[Theorem 41 of \cite{deBrecht-QuasiPolish-2013}]\label{thm: quasi-Polish space cts open surj image of cN} Suppose $Y$ is a nonempty $T_0$ space. Then $Y$ is quasi-Polish if and only if there exists a continuous open surjection $f: \cN \to Y$.
\end{theorem}

In particular, this also means that an arbitrary nonempty quasi-Polish space $Y$ is necessarily endowed with a quotient topology coming from $\cN$.

It is a consequence of theorems~\ref{thm: Polish space cts open surj image of cN}, \ref{thm: perfect Polish space cts bij image of cN}, and \ref{thm: quasi-Polish space cts open surj image of cN} that as long as a property of Polish, perfect Polish, or quasi-Polish spaces is preserved by continuous maps in the above manner, it suffices to study that property in the context of $\cN$.

Yet another consequence of these is described as follows. Given a topological space $X$, we define the collection of \textit{Borel sets} of $X$, denoted $\cB(X)$, as the smallest collection of subsets of $X$ that includes the open sets and is closed under taking complements and countable unions.\footnote{For a Polish space $X$, one can inductively define a notion of complexity of the Borel sets of $X$, called the \textit{Borel hierarchy}. (See \cite[Ch.11B]{Kechris-classical-DST-1995} for details.) This is a notion of definability used to understand the complexity of these sets in the case of the (ordinary or set-theorist's) reals.} Two spaces $(X,\cB(X))$ and $(Y,\cB(Y))$ are said to be \textit{Borel isomorphic} if, roughly, they have the same Borel structure.\footnote{See \cite[\S1.3]{Gao-2009} or \cite[Ch.12A]{Kechris-classical-DST-1995} for details.} Using Theorem~\ref{thm: Polish space cts open surj image of cN} (along with the existence of an injective Baire class-1 map from any uncountable Polish space $Y$ to $\cN$ \cite[Theorem 1.3.7]{Gao-2009}), one obtains:

\begin{corollary}[Corollary 1.3.8 in \cite{Gao-2009}] If $X$ and $Y$ are uncountable Polish spaces, then there is a Borel isomorphism from $X$ to $Y$.
\end{corollary}

By virtue of this corollary, we may define \textit{standard Borel spaces} as those whose Borel structure comes from \textit{some} uncountable Polish spaces.\footnote{Still, $\cN$ ends up being a convenient choice: See, for example, ``Since we do not lose any generality by focusing on one Polish space, we will choose to focus on the Baire space $\omega^\omega$, where this coding is easy to describe'' \cite[p.18]{Gao-2009}, as well as proofs of theorems involving standard Borel spaces throughout \cite{Gao-2009}.} For instance, quasi-Polish spaces are standard Borel (see \cite[Proposition 3.1]{Chen-QuasiPolish-2018}). In this manner, one may also define \textit{smooth Borel equivalence relations} to be those equivalence relations $E$ that are \textit{Borel reducible} to the identity relation on \textit{some} uncountable Polish space $Y$, denoted $E \leq_B \id(Y)$.\footnote{See \cite{Thomas-Schneider-CBER-2012} and \cite[\S II]{Gao-2009} for details. See also \cite{Hanson-any-function-measurable-2025}, where Hanson mentions that this is ``part of the reason why so much descriptive set theory literature is phrased explicitly in terms of Baire space, $\bN^\bN$, rather than $\bR$, despite the historical origins of the field.''}

\section{Summing Up, and a Methodological Observation}
In short, the descriptive set theorist's practice of calling the elements of $\cN$ ``reals'' is neither an abuse of language nor a mere convenience. It reflects two complementary facts that together make $\cN$ the natural exemplar of the theory. On the one hand, the very features that would disqualify $\cN$ as a stand-in for $\bR$ elsewhere in mathematics are virtues in the descriptive set theoretic context.\footnote{It is worth noting that a similar situation arises with the Cantor space $\mc{C}$: the very features that would
disqualify $\mc{C}$ as a stand-in for $\bR$ elsewhere in mathematics are virtues in the computability-theoretic context. For example, $\mc{C}$ is particularly well-suited to the study of $\Pi^0_1$-classes of reals, i.e., effectively closed subsets of $\mc{C}$ \cite[Chapter~2.19]{Downey-Hirschfeldt-2010}.} Its total disconnectedness and the homeomorphism $\cN\cong\cN^k$ render dimension invisible, so that the projective hierarchy, whose engine is projection, runs on $\cN$ without the corrections that $\bR$ (or even Cantor space $\mc{C}$) would demand. On the other hand, $\cN$ loses nothing: since every nonempty Polish space is a continuous open image of $\cN$, every perfect Polish space is a continuous bijective image, and the quasi-Polish spaces are precisely its continuous open images. Hence a property preserved under such maps may be established once, for $\cN$, and inherited by the rest. $\cN$ earns the name ``reals'' in a descriptive setting through being tractable where $\bR$ is awkward, and universal despite seeming particular. The scattered fragments of folklore can be woven into a single coherent tale.

It's worth pausing to make a methodological aside. A neglected topic in the study of mathematical technique is the importance of examples: If they are properly chosen, they can facilitate understanding and discovery.\footnote{See, for example, \cite{Antos-2024} where the authors study explanatory proofs in mathematics through an extended case study from DST. Antos is currently carrying out more detailed case studies of examples in mathematics in work in progress.} The role of Baire space in the folklore of DST provides a compelling example: Even in cases where the aim is to demonstrate results that hold for every perfect Polish space, there is one specific instance which incorporates the necessary information in an especially fruitful and illuminating way.

\bibliography{References}
\bibliographystyle{alpha}

\bigskip
\noindent\textsc{Dhruv Kulshreshtha}\\
Department of Mathematics, University of Wisconsin--Madison\\
\textit{Email:} \href{mailto:dkulshreshth@wisc.edu}{\texttt{dkulshreshth@wisc.edu}}

\medskip
\noindent\textsc{Jamie Tappenden}\\
Department of Philosophy, University of Michigan, Ann Arbor\\
\textit{Email:} \href{mailto:tappen@umich.edu}{\texttt{tappen@umich.edu}}

\end{document}